%February 11, 2004
\input amstex
\documentstyle{amsppt}
\magnification=1200
\hoffset=-0.5pc
\vsize=57.2truepc
\hsize=37truepc
\spaceskip=.5em plus.25em minus.20em
\define\Bobb{\Bbb}
\define\rrtimes{\rtimes}
\define\Ho{{\roman H}}
\define\differen {1}
\define\modpcoho{2}
\define\huebkade{3}

\nologo
\topmatter
\title
On the cohomology of the holomorph of
a finite cyclic group \endtitle
\author
Johannes Huebschmann 
\endauthor
\date
{February 11, 2004}
\enddate
\abstract
{We determine the mod 2 cohomology algebra of the holomorph of any finite 
cyclic group whose order is a power of 2.}
\endabstract
\affil 
Universit\'e des Sciences et Technologies de Lille
\\
UFR de Math\'ematiques
\\
CNRS-UMR 8524
\\
F-59 655 VILLENEUVE D'ASCQ C\'edex/France
\\
Johannes.Huebschmann\@math.univ-lille1.fr
\endaffil
\address{\noindent
USTL, UFR de Math\'ematiques, 
F-59 655 VILLENEUVE D'ASCQ C\'edex/France
\newline\noindent
Johannes.Huebschmann\@math.univ-lille1.fr
}
\endaddress
\endkeywords
\endtopmatter
\document
\rightheadtext{Cohomology of the holomorph}

\beginsection 1. Outline

The {\it holomorph\/} of a group is the semi-direct product of the group with
its automorphism group,
with respect to the obvious action.
The automorphism group of a
non-trivial finite cyclic group of order $r$ is well known
to be cyclic
if and only if the number $r$ is of the kind
$r=4, r=p^{\rho}, r=2p^{\rho}$ where
$p$ is an odd prime;
in these cases, the holomorph is thus a split metacyclic group.
The mod $p$ cohomology algebra of an arbitrary metacyclic group
has been determined in \cite\modpcoho.
In this note we will 
determine the mod 2 cohomology of the 
holomorph of a cyclic group whose order is a power of 2.
Since for odd $p$ the $p$-primary part of the
automorphism group of any cyclic group is cyclic,
we thus get a complete description of the
mod $p$ cohomology of the holomorph of any (finite) cyclic group.
I am indebted to Fred Cohen for having asked whether 
the cohomology of the holomorph of a cyclic group whose order is a power of 2,
this holomorph not being metacyclic when the order of the cyclic group
is at least 8,
may be determined by means of the methods I developed in the 80's.

Let $\rho \geq 2$ and let $N$ be the cyclic group of order $r=2^{\rho}$.
As usual, we will identify the automorphism group of $N$ with the group
$(\Bobb Z\big /2^{\rho})^*$ of units of $\Bobb Z\big /2^{\rho}$,
the latter being viewed as a commutative ring.
Consider the holomorph 
$$
G=(\Bobb Z\big /2^{\rho}) \rrtimes (\Bobb Z\big /2^{\rho})^*
\tag1.1
$$
of $N$.
For $\rho = 2$, this group  comes down to the dihedral group,
the group $(\Bobb Z\big /4)^*$ being cyclic of order 2,
generated by the class of $-1$.
Henceforth we suppose that
$\rho \geq 3$.
Now the group $(\Bobb Z\big /2^{\rho})^*$
decomposes as a direct product of 
a copy of $\Bobb Z\big /2$, 
generated by the class of $-1$, and
a copy of $\Bobb Z\big /2^{\rho-2}$,
generated by the class of $5$.
Write $s =2^{\rho-2}$.
The cyclic groups being written multiplicatively,
the semi-direct product (1.1) has thus the presentation
$$
\langle x,y,z;\ y^r = 1,\ x^s = 1,
\ xyx^{-1} = y^5,
\ zyz^{-1} = y^{-1},
\ [x,z] =1,
\ z^2=1
 \rangle,
\tag1.2
$$
the normal cyclic subgroup
$N$ being generated by $y$.

Denote by
$K_x$ and $K_z$
the cyclic subgroups
of order $s$ and $2$ generated by $x$
and $z$, respectively.
The mod 2 cohomology algebra $\Ho^*(N,\Bobb Z /2)$
is well known to be generated by a class
${\omega_y \in \Ho^1(N,\Bobb Z /2)}$
and a class
${c_y \in \Ho^2(N,\Bobb Z /2)}$,
subject to the relation $\omega^2_y = \frac r2 c_y$.
Likewise,
the mod 2 cohomology algebra $\Ho^*(K_z,\Bobb Z /2)$
is freely generated by a class
${\omega_z \in \Ho^1(K_x,\Bobb Z /2)}$,
and
the mod 2 cohomology algebra $\Ho^*(K_x,\Bobb Z /2)$
is generated by certain classes
${\omega_x \in \Ho^1(K_x,\Bobb Z /2)}$ and ${c_x \in \Ho^2(K_x,\Bobb Z /2)}$,
subject to the relation $\omega^2_x = \frac s2 c_x$.

While the cohomology of $G$ may be determined
by means of the obvious split
extension of $N$ by 
$K_x \times K_z$, a more economical approach
to solving the resulting multiplicative extension
problem involves 
the group 
$\widehat G$  
which is given by the presentation
$$
\langle x,y,z;\ y^r = 1,
\ xyx^{-1} = y^5,
\ zyz^{-1} = y^{-1},
\ [x,z] =1,
\ z^2=1
 \rangle
\tag1.3
$$
and projects onto $G$ in an obvious fashion;
here the term \lq\lq more economical\rq\rq\ will be justified
in Remark 3.7 below.
Let $\omega_1$  be the homomorphism from $G$ 
to $\Bobb Z/2$ which sends $y$ to the generator
of $\Bobb Z/2$ and is trivial on the two other generators $x$ and $z$ 
of $G$. Likewise, 
abusing the notation $\omega_x$  and $\omega_z$ somewhat, let
$\omega_x$  and $\omega_z$
be the homomorphisms from $G$ 
to $\Bobb Z/2$ which send $x$ and
$z$, respectively,
to the generator 
of $\Bobb Z/2$ and are trivial on the respective two other generators.
These homomorphisms are well defined and
plainly factor through 
the homomorphisms on $K_x$ and $K_z$ 
which yield the cohomology classes of these groups
denoted by the same symbols.
Abusing this notation further, we will denote the composites
of $\omega_1$,  $\omega_x$,  and $\omega_z$
with the projection from
$\widehat G$ to $G$ by the same symbols as well.
Consider the group $\widetilde G$ given by the presentation
$$
\langle x,\widetilde y,z;\ \widetilde y^{2r} = 1,
\ x\widetilde y x^{-1} = \widetilde y^5,
\ z\widetilde yz^{-1} = \widetilde y^{-1},
\ [x,z] =1,
\ z^2=1
 \rangle.
$$
The obvious homomorphism from $\widetilde G$ to $\widehat G$
which sends $\widetilde y$ to $y$ and
the other generators to the generators denoted by the same symbols yields
a central extension
$$
0
@>>>
\Bobb Z/2
@>>>
\widetilde G
@>>>
\widehat G
@>>>
1
$$
the class of which we denote by $c_2 \in \Ho ^2(\widehat G,\Bobb Z/2)$.
Since the extension, restricted to $N = \Bobb Z/r$, yields the extension
representing the class $c_y \in \Ho ^2(N,\Bobb Z/2)$,
the class $c_2$ restricts to $c_y \in \Ho ^2(N,\Bobb Z/2)$;
likewise, $\omega_1$
restricts to $\omega_y \in \Ho^1(N,\Bobb Z/2)$.

\proclaim{Proposition 1.4}
As a graded commutative algebra,
${\Ho ^*(\widehat G,\Bobb Z /2)}$ is generated by
$\omega_x$, $\omega_z$,  $\omega_1$, $c_2$, subject to the relations
$$
\align
\omega_x^2 &= 0
\tag1.4.1
\\
\omega_1^2&= \omega_z \omega_1.
\tag1.4.2
\endalign
$$
\endproclaim

We note that, in characteristic 2, there is
no difference between graded commutative and commutative. 
Here is our main result.

\proclaim{Theorem 1.5}
The mod {\rm 2} cohomology of the group $G$ has classes
$$
\omega_3\in \Ho^3(G,\Bobb Z /2),
\ c_4\in \Ho^4(G,\Bobb Z /2)
$$
which, under inflation, go to the classes
$$
(\omega_1+\omega_z) c_2\in \Ho^3(\widehat G,\Bobb Z /2),
\ c_2^2\in \Ho^4(\widehat G,\Bobb Z /2)
$$
such that the mod {\rm 2} cohomology algebra of $G$
is generated by
$$
c_x,\ \omega_x,\ \omega_z,\ \omega_1,\  \omega_3,\ c_4,
$$
subject to the relations 
$$
\align
\omega_1 \omega_3 &= 0\tag1.5.1\\
c_x \omega_1&= 0\tag1.5.2\\
\omega_x^2&= \tfrac s2 c_x
\tag1.5.3\\
\omega_1^2&=\omega_z \omega_1
\tag1.5.4\\
\omega_3^2 &=\omega_z\omega_1 c_4 + \omega_z^2c_4 .
\tag1.5.5
\endalign
$$ 
Requiring that
the classes $\omega_3$ and $c_4$ 
restrict to zero in the cohomology of the abelian subgroup of $G$ generated
by $x$ and $z$ determines these classes uniquely.
\endproclaim

In Section 3 of \cite\modpcoho,
by means of {\it homological perturbation theory\/}, 
we constructed  a free resolution
for an arbitrary metacyclic group
from a presentation thereof;
see  
\cite\modpcoho\ and \cite\huebkade\ 
for comments about 
homological perturbation theory
and for more references.
The construction 
may be extended to 
that of a free resolution
for the group $G$ from the presentation (1.2),
and an explicit description of the 
cohomology as well as of various
spectral sequences
may be
derived from it.
In this paper we give an alternate somewhat simpler approach
which essentially reduces 
the requisite structural insight
to insight into the resolutions for various metacyclic groups
which we constructed in \cite\modpcoho.
In particular, we shall detect the behaviour of various differentials
in certain spectral sequences by inflation to spectral sequences
of related group extensions and we shall determine the multiplicative
relations in a similar fashion.

\medskip\noindent {\bf 2. The additive structure}
\smallskip
\noindent
(2.1) The subgroup $G_x$ of $G$ generated by $x$ and $y$ is metacyclic
with presentation 
$$
\langle x,y;\quad y^r = 1,\quad x^s = 1,
\quad xyx^{-1} = y^5 \rangle
\tag2.1.1
$$
and fits into the split extension
$$
\roman e_x \colon
1
@>>>
N
@>>>
G_x
@>>>
K_x
@>>>
1.
\tag2.1.2
$$
The number 
$$
\frac{5^{2^{\rho-2}}-1}{2^{\rho}}=
\frac{(5^{2^{\rho-3}})^2-1}{2^{\rho}}=
\frac{5^{2^{\rho-3}}-1}{2^{\rho-1}}
\frac{5^{2^{\rho-3}}+1}2
$$
is odd, and 
we are in the situation
of Theorem C of
\cite\modpcoho\ (with the notation of that
Theorem: for the case $t=5$ and $p=2$).
In particular,
the mod 2 cohomology spectral sequence 
$(\roman E_r^{s,t}(\roman e_x), d_r)$
for the extension $\roman e_x$
has $d_1$ equal to zero and,
as a commutative algebra,
$\roman E_2(\roman e_x)$
is generated by the classes
$\omega_x$,  $c_x$, $\omega_y$, $c_y$,
subject to the relations 
$$
\omega^2_x = \frac s2 c_x, \  \omega^2_y = \frac r2 c_y.
$$
In view of Lemma 4.2.1 of 
\cite\modpcoho,
in this spectral sequence, for $j \geq 1$,
$$
d_2(c_y^j) = j c_x c^{j-1}_y \omega_y,
\ 
d_2(c_y^j \omega_y) = 0
\tag2.1.3
$$
and, as explained on pp. 84 ff. of \cite\modpcoho,
the spectral sequence collapses from $\roman E_3(\roman e_x)$.
For later reference we recall
the following simplified version of Theorem C of
\cite\modpcoho.

\proclaim{Proposition 2.1.4}
The cohomology spectral sequence of the group extension
$\roman e_x$
collapses from $\roman E_3(\roman e_x)$, and 
$\Ho ^*(G_x,\Bobb Z /2)$ has classes
$\omega_3 \in \Ho^3(G_x,\Bobb Z /2)$ 
and
$c_4 \in \Ho^4(G_x,\Bobb Z /2)$
which restrict to the classes
${\omega_y c_y}$ and ${c_y^2}$, respectively,
in ${\Ho^*(N,\Bobb Z /2)}$
so that, as a commutative algebra,
${\Ho ^*(G_x,\Bobb Z /2)}$ is generated by
$c_x$, $\omega_x$, $\omega_1$,  $\omega_3$, $c_4$,
subject to the following relations:
$$
\align
\omega_1 \omega_3 &= 0\tag2.1.4.1\\
 c_x \omega_1&=0\tag2.1.4.2\\
\omega_x^2&= \tfrac s2 c_x
\tag2.1.4.3\\
\omega_1^2&=0           
\tag2.1.4.4\\
\omega_3^2 &=0 .
\tag2.1.4.5
\endalign
$$ 
\endproclaim

\noindent
(2.1.5) Later we shall have to resolve various ambiguities related with the 
multiplicative structure of the group $G$ which we are 
really interested in. We now describe the corresponding
ambiguities for the group $G_x$.
We suppose that a choice of the classes 
$\omega_3 \in \Ho ^3(G_x,\Bobb Z /2)$ and $c_4 \in \Ho ^4(G_x,\Bobb Z /2)$
has been made. In dimensions 1,2,3,4 the following monomials then constitute
bases:
$$
\align
\Ho^1(G_x,\Bobb Z /2) \colon &\quad \omega_1, \omega_x
\\
\Ho^2(G_x,\Bobb Z /2) \colon &\quad \omega_1\omega_x, c_x
\\
\Ho^3(G_x,\Bobb Z /2) \colon &\quad \omega_x c_x, \omega_3
\\
\Ho^4(G_x,\Bobb Z /2) \colon &\quad c^2_x, c_4,
\omega_3 \omega_x
\endalign
$$
Hence requiring that, under the restriction to $\roman H^*(K_x,\Bobb Z/2)$,
$\omega_3$ and $c_4$ go to zero determines $\omega_3$ uniquely and
$c_4$ up to the term $\omega_3 \omega_x$.

\noindent
(2.2) The subgroup $G_z$ of $G$ generated by
$z$ and $y$ is metacyclic of the form
$G_z=N \rrtimes K_z$
and thus fits into the split extension
$$
\roman e_z \colon
1
@>>>
N
@>>>
G_z
@>>>
K_z
@>>>
1.
$$
The (images of the) classes $\omega_z$, $\omega_1$, $c_2$
in $\Ho ^*(G_z,\Bobb Z /2)$
generate the  constituent $\roman E_2(\roman e_z)$ of
the mod 2 cohomology spectral sequence of $\roman e_z$ 
which necessarily collapses from $\roman E_2$
whence, as a (graded) commutative algebra,
$\Ho ^*(G_z,\Bobb Z /2)$ 
is generated by these classes.
By Theorem B of \cite\modpcoho, 
the relation (1.4.2) is defining, 
that is, ${\Ho ^*(G_z,\Bobb Z /2)}$ 
is generated by these classes
subject to the relation (1.4.2).

\noindent
(2.3) The projection from $\widehat G$ to the cyclic group generated by $x$
yields the group extension
$$
\widehat{\roman e}\colon
1
@>>>
G_z
@>>>
\widehat G
@>>>
\Bobb Z
@>>>
1.
\tag2.3.1
$$
We now give a somewhat more precise version of Proposition 1.4
in the introduction.

\proclaim{Proposition 2.3.2}
The spectral sequence of the extension $\widehat {\roman e}$ 
collapses from $\roman E_2$ and,
as a graded commutative algebra,
${\Ho ^*(\widehat G,\Bobb Z /2)}$ is generated by
$\omega_x,\ \omega_z, \ \omega_1,\ c_2$,
subject to the relations
{\rm (1.4.1)} and {\rm (1.4.2)}.
\endproclaim

\demo{Proof}
Since every multiplicative generator of the cohomology of $G_z$
is in the image of the restriction map from $\widehat G$ to $G_z$,
the multiplicative properties
of the spectral sequence of the extension  
$\widehat {\roman e}$
entail that the spectral sequence collapses from $\roman E_2$.
The relation (1.4.1) holds for degree reasons.
Moreover, since in $\roman H^*(G_z,\Bobb Z/2)$ the relation
(1.4.2) holds,
for suitable coefficients $a,b \in \Bobb Z/2$,
$$
\omega_1^2= \omega_z \omega_1 + a \omega_z\omega_x +b \omega_1\omega_x .
$$
Restricting this relation to the subgroup 
of $\widehat G$ (of the kind $K_z \times \Bobb Z$)
generated by
$x$ and $z$ we see that $a=0$.
Likewise,
restricting the relation 
to the subgroup 
$\widehat G_x$ 
of $\widehat G$ 
(of the kind $N \rrtimes \Bobb Z$) generated by
$x$ and $y$, by virtue of Proposition 2.1.4,
since $\omega_1^2 = 0$ in 
$\roman H^2(\widehat G_x,\Bobb Z/2)$,
we conclude that $b=0$. \qed
\enddemo

\noindent
(2.4) The projection from $G$ to the cyclic group generated by $x$
yields the group extension
$$
\roman e
\colon
1
@>>>
G_z
@>>>
G
@>>>
K_x@>>>
1
\tag2.4.1
$$
in an obvious fashion.
In view of (2.2) above,
as a graded $\Bobb Z/2$-algebra, the cohomology
$\roman H^*(G_z,\Bobb Z/2)$ 
of $G_z$ is generated by
$\omega_z$, $\omega_1$, $c_2$,
subject to the relation (1.4.2).

\proclaim{Proposition 2.4.2}
The cohomology spectral sequence 
$(\roman E_r(\roman e), d_r)$ of the group extension $\roman e$
has 
$$
d_2(\omega_1)=0, \quad d_2(c_2) =   \omega_1 c_x
$$
and collapses from $\roman E_3(\roman e)$. Furthermore,
${\Ho ^*(G,\Bobb Z /2)}$ has a class $\omega_3 \in
\roman H^3(G,\Bobb Z /2)$
which restricts to 
$(\omega_1 + \omega_z) c_2 \in \roman H^3(G_z,\Bobb Z /2)$
and a class
$c_4 \in \roman H^4(G,\Bobb Z /2)$
which restricts to $c_2^2 \in \roman H^4(G_z,\Bobb Z /2)$
so that,
as a
commutative algebra,
${\Ho ^*(G,\Bobb Z /2)}$ is generated by
$c_x,\ \omega_x,\ \omega_z,\ \omega_1,\  \omega_3,\ c_4$.
\endproclaim

\demo{Proof}
The class $\omega_1$ is manifestly an infinite cycle.
For suitable coefficients $\eta,\epsilon \in \Bobb Z/2$,
the spectral sequence $(\roman E_r(\roman e), d_r)$ necessarily has
$$
d_2(c_2) =  \eta \omega_1 c_x + \epsilon \omega_z  c_x.
$$
Consider the morphism
$$
\CD
1
@>>>
N
@>>>
G_x
@>>>
K_x
@>>>
1
\\
@.
@VVV
@VVV
@VV{\roman{Id}}V
@.
\\
1
@>>>
G_z
@>>>
G
@>>>
K_x
@>>>
1
\endCD
\tag2.4.3
$$
of extensions from $\roman e_x$ to $\roman e$.
Restricting the spectral sequence $(\roman E_r(\roman e), d_r)$ 
to the extension $\roman e_x$ we deduce from (2.1.3) that $\eta=1$.
Likewise, restricting the spectral sequence to the extension
$$
\roman e^{\times}\colon
1
@>>>
K_z
@>>>
K_z \times K_x
@>>>
K_x@>>>
1
\tag2.4.4
$$
we see that $\epsilon = 0$. Consequently,
in the spectral sequence $(\roman E_r(\roman e), d_r)$, 
$$
d_2(c_2) = c_x \omega_1 
\tag2.4.5
$$
whence, in view of the multiplicative structure of the spectral sequence,
for $j \geq 1$,
$$
d_2(c^j_2) = j c^{j-1}_2c_x \omega_1,
\ 
d_2(c^j_2 \omega_1) = j c^{j-1}_2c_x \omega^2_1 
=j c^{j-1}_2c_x \omega_1 \omega_z .
\tag2.4.6
$$
In particular, $d_2c_2^2 = 0$. 
We assert that $c_2^2$ is an infinite cycle.
Indeed, the value of
$d_3c_2^2$ lies in $\roman E_3^{3,2}(\roman e)$.
Restricting the spectral sequence $(\roman E_r(\roman e), d_r)$
to the two extensions $\roman e_x$ and $\roman e^{\times}$ 
we deduce that $d_3c_2^2$ is zero, and we are left with
$$
d_4c_2^2 \in \roman E_4^{4,1}(\roman e) .
$$
The same kind of detection principle as that used above shows that
the only possibility for
$d_4c_2^2$ to be non-zero is that, for 
a suitable coefficient $a \in \Bobb Z/2$,
$$
d_4c_2^2 =a \omega_z^3 \omega_x \omega_y.
$$
However, the obvious projection map 
yields a morphism of extensions from
$\widehat {\roman e}$ to 
$\roman e$
and hence a morphism of spectral sequences
via the inflation map.
By Proposition 2.3.2, the mod 2 cohomology spectral sequence 
$(\roman E_r^{s,t}(\widehat {\roman e}), d_r)$
collapses from $\roman E_2$.
Hence the coefficient $a$ is zero,
that is, $d_4c_2^2$ is zero.
Finally, since the extension $\roman e$ is split,
a non-zero differential
of the spectral sequence
$(\roman E_r^{s,t}(\roman e), d_r)$
cannot hit the baseline.
Consequently
$d_5c_2^2$ is zero whence
$c_2^2$ is an infinite cycle.

Likewise 
$(\omega_1 + \omega_z) c_2$
is an infinite cycle.
The multiplicative properties of the spectral sequence 
$(\roman E_r^{s,t}(\roman e), d_r)$
imply that this spectral sequence collapses from $\roman E_3(\roman e)$.
We conclude that, apart from the multiplicative generators
$c_x$, $\omega_x$, $\omega_z$, $\omega_1$,
the cohomology algebra
$\roman H^*(G,\Bobb Z /2)$
has  an additional generator 
$\omega_3 \in
\roman H^3(G,\Bobb Z /2)$
which restricts to 
$(\omega_1 + \omega_z) c_2 \in \roman H^3(G_z,\Bobb Z /2)$
and another generator
$c_4 \in \roman H^4(G,\Bobb Z /2)$
which restricts to $c_2^2 \in \roman H^4(G_z,\Bobb Z /2)$
such that, as an algebra,
$\roman H^*(G,\Bobb Z /2)$
is generated by these six generators. \qed
\enddemo

\noindent
(2.4.7) 
We suppose that a choice of the classes 
$\omega_3 \in \Ho ^3(G,\Bobb Z /2)$ and $c_4 \in \Ho ^4(G,\Bobb Z /2)$
has been made. In dimensions 3 and 4 the following monomials then constitute
bases:
$$
\align
\Ho^3(G,\Bobb Z /2) \colon &\quad 
\omega_x c_x, \omega_3;
\omega_z\omega_1\omega_x, \omega_z c_x;
\omega^2_z\omega_1, \omega^2_z\omega_x, \omega_z^3
\\
\Ho^4(G,\Bobb Z /2) \colon &\quad 
c^2_x, c_4, \omega_3 \omega_x;
\omega_z\omega_x c_x, \omega_z\omega_3;
\omega^2_z\omega_1\omega_x, \omega^2_z c_x;
\omega^3_z\omega_1, \omega^3_z\omega_x, \omega_z^4
\endalign
$$
Requiring that $\omega_3$ and $c_4$
restrict to zero in the cohomology of the abelian subgroup
of $G$ generated by $x$ and $z$
leaves the following ambiguities:
$$
\align
\omega_3:&\quad \omega^2_z\omega_1, \omega_z\omega_1\omega_x
\\
c_4:&\quad \omega^3_z\omega_1,
\omega^2_z\omega_1 \omega_x,
\omega_3 \omega_x .
\endalign
$$
Requiring that, furthermore,
the classes $\omega_3$ and $c_4$ go,
under inflation, to the classes
$$
(\omega_1+\omega_z) c_2\in \Ho^3(\widehat G,\Bobb Z /2),
\ c_2^2\in \Ho^4(\widehat G,\Bobb Z /2),
$$
respectively, as asserted in Theorem 1.5,
determines the classes $\omega_3 
\in \Ho^3(G,\Bobb Z /2)$ and
${c_4\in \Ho^4(G,\Bobb Z /2)}$ uniquely,
since the classes
$\omega^2_z\omega_1, \omega_z\omega_1\omega_x$
and
$\omega^3_z\omega_1,\omega^2_z\omega_1 \omega_x,\omega_3 \omega_x$
go to linearly independent classes in
$\Ho^*(\widehat G,\Bobb Z /2)$.
Henceforth we suppose that these choices
for $\omega_3$ and $c_4$ which determine them uniquely have been made.
Then, under {\rm (2.4.3)},
the generator $\omega_3 \in\roman H^3(G,\Bobb Z /2)$ 
goes to the generator
in $\roman H^3(G_x,\Bobb Z /2)$
denoted in
{\rm (2.1.4)} by the same symbol
and $c_4\in\roman H^4(G,\Bobb Z /2)$
goes to one of the two possible choices for
the class in $\roman H^4(G_x,\Bobb Z /2)$
denoted $c_4$ in {\rm (2.1.4)} 
and hence resolves the ambiguity left open in (2.1.5).

These observations yield, in particular, a proof of Theorem 1.5, 
as far as the additive structure is concerned.

\noindent
{\smc Remark 2.4.8.} The value 
$d_2(c_2) =   \omega_1 c_x$ of the class $c_2$
under the differential $d_2$ in the spectral sequence
of the extension $\roman e$ may also be determined by means of
Theorem 1 in \cite\differen.

\medskip\noindent{\bf 3. The multiplicative structure}
\smallskip\noindent
The obvious projection from $\widehat G$ 
to $G$ yields a central extension
$$
1
@>>>
\Bobb Z
@>>>
\widehat G
@>>>
G@>>>
1
\tag3.1
$$
in an obvious fashion.
We realize this central extension geometrically by a fiber bundle of the kind
$$
B \widehat G
@>>>
B G
@>>>
BS^1
\tag3.2
$$
where
the map from $B G$ to
$BS^1$ represents the cohomology class $c_x$;
here $S^1$ is the circle group.
By Proposition 2.3.2, as a graded $\Bobb Z/2$-algebra, the cohomology
$\roman H^*(\widehat G,\Bobb Z/2)$ of the fiber  
is generated by
$\omega_x$, $\omega_z$, $\omega_1$, $c_2$,
subject to the relations $\omega_1^2=\omega_z \omega_1$
and $\omega_x^2=0$.

Consider the mod 2 spectral sequence $(\roman E_r^{s,t}, d_r)$ for (3.2). 
It has
$$
\roman E_2^{*,t} = \Bobb Z/2 [c_x] \otimes 
\roman H^t(\widehat G,\Bobb Z/2) 
$$
and, as usual, we will write
$\roman E_{\infty} = \roman E_{\infty}(\roman H^*(G,\Bobb Z/2))$
for the associated graded algebra
coming from the filtration 
in terms of powers of the multiplicative generator $c_x$.
Given $p,q \geq 0$, this filtration has the form
$$
\roman H^{2p,q} \subseteq \roman H^{2p-2,q+2} \subseteq \dots 
\subseteq
\roman H^{0,q+2p} =\roman H^{q+2p}(G,\Bobb Z/2).
$$
Since $\omega_z$ and $\omega_1$
are restrictions of cohomology classes of $G$,
in the spectral sequence, they are infinite cycles.
Moreover, in view of (2.4.5), the spectral sequence has
$$
d_2(c_2) =   \omega_1 c_x. 
\tag3.3
$$
More precisely, for degree reasons, 
for suitable coefficients $a,b,c$ in $\Bobb Z/2$,
$$
d_2(c_2) = a  \omega_1 c_x + b \omega_x  c_x + c \omega_z  c_x.
$$
Restricting the extension (3.1) to the extension
$$
1
@>>>
\Bobb Z
@>>>
\widehat G_x
@>>>
G_x@>>>
1
$$
we see that $a=1$ and $b=0$ and restricting 
the extension (3.1) to the trivial extension
$$
1
@>>>
1
@>>>
G_z
@>>>
G_z@>>>
1
$$
so that $B S^1$ is replaced with a point we see that $c=0$.
Likewise
$$
d_2(\omega_1c_2) =\omega^2_1 c_x = \omega_z \omega_1 c_x,
\quad
d_2(\omega_zc_2) =\omega_z \omega_1 c_x,
$$
and $(\omega_1 + \omega_z) c_2$
is an infinite cycle, indeed, arises from the class
$$
\omega_3 \in \roman H^3(G,\Bobb Z/2)= \roman H^{0,3}
$$ 
under the surjection from $\roman H^{0,3}$ to $\roman E_{\infty}^{0,3}$.
Furthermore, 
$d_2(c^2_2) =0$
whence $c^2_2$ is an infinite cycle.
It arises from the class $c_4 \in \roman H^4(G,\Bobb Z/2)=\roman H^{0,4}$
under the surjection from $\roman H^{0,4}$ to $\roman E_{\infty}^{0,4}$.
Abusing notation somewhat, we will write
$$
\omega_3=(\omega_1 + \omega_z) c_2,\quad 
c_4 =c^2_2.
$$

\proclaim{Proposition 3.4}
The associated graded algebra
$\roman E_{\infty}(\roman H^*(G,\Bobb Z/2))$
is generated by
$c_x$, $\omega_x$, $\omega_z$, $\omega_1$,  $\omega_3$, $c_4$,
subject to the relations 
$$
\align
\omega_1 \omega_3 &= 0 
\tag1.5.1$(\infty)$\\
c_x \omega_1&= 0\tag1.5.2$(\infty)$\\
\omega_x^2&= \tfrac s2 c_x
\tag1.5.3\\
\omega_1^2&=\omega_z \omega_1  
\tag1.5.4$(\infty)$\\
\omega_3^2 &=\omega_z\omega_1 c_4 + \omega_z^2c_4 .
\tag1.5.5$(\infty)$
\endalign
$$ 
\endproclaim

\demo{Proof}
Indeed, since $\omega^2_1 = \omega_1 \omega_z \in 
\roman H^2(\widehat G,\Bobb Z/2)$, we have
$$
\omega_1\omega_3 = \omega_1(\omega_1 + \omega_z) c_2 =
(\omega_1^2 + \omega_1\omega_z) c_2 =0,
$$
that is, the relation (1.5.1$(\infty)$) holds.
The relation (1.5.2$(\infty)$) 
is an immediate consequence of (3.3).
The relation (1.5.3) holds already in $\roman H^2(K_x,\Bobb Z/2)$,
and the relation (1.5.4$(\infty)$) holds in 
$\roman H^2(\widehat G,\Bobb Z/2)$.
Finally,
$$
\omega_3^2 = \omega^2_1 c_4+ \omega^2_z c_4
=\omega_1\omega_z c_4+ \omega^2_z c_4,
$$
whence the relation (1.5.5$(\infty)$) is satisfied. \qed 
\enddemo

We now complete the proof of Theorem 1.5 by
verifying that the relations 
(1.5.1$(\infty)$), (1.5.2$(\infty)$), (1.5.4$(\infty)$),
and (1.5.5$(\infty)$) hold in $\roman H^*(G,\Bobb Z/2)$.
This amounts to showing that the multiplicative extension problem 
with reference to the filtration 
coming from powers of the generator $c_x$ is trivial.
To this end  we examine the possible ambiguities
of the corresponding relations
in $\roman H^*(G,\Bobb Z/2)$.

The relation (1.5.4$(\infty)$) arises
from a relation of the kind
$$
\omega^2_1 = \omega_z \omega_1 + a c_x
$$
for some $a \in \Bobb Z/2$. Restricting this relation
to the abelian subgroup of $G$ generated by $x$ and $z$ we find that $a=0$
whence the relation (1.5.4) holds.

Likewise the relation (1.5.2$(\infty)$) arises from
a relation of the kind
$$
c_x \omega_1 = a \omega_z c_x +b \omega_x c_x,
$$
for suitable coefficients $a, b \in \Bobb Z/2$.
When we restrict this relation  to the abelian subgroup 
of $G$ generated by $x$ and $z$ we find $a=b=0$,
that is, the relation (1.5.2) holds.

Since $\omega_1 \omega_3$ restricts to zero in the cohomology of
$\widehat G$,
the relation (1.5.1$(\infty)$) arises from
a relation of the kind
$$
\omega_1 \omega_3 = 
b(\omega_x,\omega_z,\omega_1)c_x + d c_x^2
\tag1.5.1$'$
$$
where $b$ is a quadratic polynomial in $\omega_x,\omega_z,\omega_1$
and where $d \in \Bobb Z/2$.
The group $\roman H^2(\widehat G,\Bobb Z/2)$ has dimension 4, a basis
being given by
$$
\omega_1 \omega_z,\  \omega_1 \omega_x,\  \omega_z^2,\  \omega_z \omega_x;
\tag3.5
$$
notice that $\omega_x^2 = 0$. 
The relations $\omega^2_1 = \omega_1\omega_z$ 
and $\omega_1 c_x=0$ in the cohomology of $G$
have already been established.
Thus only the monomials $\omega_x \omega_z$ 
and $\omega_z^2$
may yield non-zero contributions
to $b(\omega_x,\omega_z,\omega_1)c_x$, that is,
for some $\beta,\gamma \in \Bobb Z/2$,
$$
\omega_1 \omega_3 = 
\beta \omega_z\omega_xc_x + \gamma \omega^2_z c_x + d c^2_x.
\tag 1.5.1$''$
$$
When we restrict the relation (1.5.1$''$) to the subgroup 
of $G$ (of the kind $K_z \times K_x$) generated
by $x$ and $z$, since $\omega_1$ restricts to zero, we see that 
$\beta = \gamma=0$ and $d=0$,
that is, the relation (1.5.1) holds.

We will now resolve the ambiguities for the relation (1.5.5$(\infty)$).
The monomials
$$
c_4, \ 
\omega_1 \omega_z^2 \omega_x, \
\omega_z^3 \omega_x, \
\omega_z^3 \omega_1, \
\omega_3 \omega_1, \
\omega_3 \omega_z, \
\omega_3 \omega_x
\tag3.6
$$
constitute a basis of $\roman H^4(\widehat G,\Bobb Z/2)$.
Since, in the cohomology of $G$,  $\omega_1 c_x = 0$,
non-zero contributions to the square $\omega_3^2$ may at first arise
from the products of each of the monomials 
$c_4$, $\omega_z^3 \omega_x$, $\omega_3 \omega_z$,
$\omega_3 \omega_x$ with $c_x$.
However, $\omega_3^2$ restricts to zero in the 
cohomology of the metacyclic subgroup
$G_x$ of $G$ generated by $x$ and $y$, and the monomials
$c_4 c_x$ 
and $\omega_3 \omega_x c_x$
are non-zero in $\roman H^6(G_x,\Bobb Z/2)$.
Likewise, $\omega_3^2$ restricts to zero in the 
cohomology of the abelian subgroup
$K_x \times K_z$ of $G$ generated by $x$ and $z$, and
$\omega_z^3 \omega_x c_x$ is non-zero in 
$\roman H^6(K_x \times K_z,\Bobb Z/2)$.
Hence, among the monomials which are {\it linear\/} in $c_x$,
the square $\omega_2^3$ may involve
at most the term $\omega_3 \omega_z c_x$.

We now show that monomials
which are {\it quadratic\/} in $c_x$
cannot yield a non-zero contribution to $\omega_3^2$.
Indeed,
these monomials arise from the monomials (3.5) since
these constitute a basis of $\roman H^2(\widehat G,\Bobb Z/2)$.
Among the four monomials of this kind, the monomials
$\omega_1 \omega_z$ and
$\omega_1 \omega_x$
will not contribute anything non-zero since
$\omega_1 c_x = 0$.
Furthermore, the monomials
$\omega_z^2c_x^2$ and $\omega_z \omega_x c_x^2$, 
restricted to the abelian subgroup
$K_x \times K_z$ of $G$ generated by $x$ and $z$,
yield non-zero elements of
$\roman H^6(K_x \times K_z,\Bobb Z/2)$.
Thus, since the square $\omega_2^3$ 
restricts to zero in $\roman H^6(K_x \times K_z,\Bobb Z/2)$,
$\omega_2^3$  cannot involve
a monomial which is {\it quadratic\/} in $c_x$.

Summing up, we conclude that,
for some $\gamma \in \Bobb Z/2$, possibly zero,
$$
\omega_3^2 =\omega_z\omega_1 c_4 + \omega_z^2c_4
+
\gamma \omega_z \omega_3 c_x.
\tag1.5.5$'$
$$
The following argument rules out the possibility  of a non-zero $\gamma$.
Consider the subgroup $A$ of $G$ generated by
$x$, $y^{\frac r2}$, and $z$.
This subgroup decomposes as a direct product
$A \cong \Bobb Z/2 \times K_z \times K_x$ where
the unlabelled copy of
$\Bobb Z/2$ is generated by 
$y^{\frac r2}$.
As a graded commutative
algebra, the mod 2 cohomology of this group is generated by
a class $\omega \colon A \to \Bobb Z/2$ which is the obvious
projection to the unlabelled copy of $\Bobb Z/2$,
together with the classes $\omega_z$, $\omega_x$, and $c_x$,
subject to the relation  $\omega^2_x=0$,
where the notation is abused.
Under the restriction map from $\roman H^*(G,\Bobb Z/2)$
to
$\roman H^*(A,\Bobb Z/2)$,
the classes $\omega_z$, $\omega_x$, and $c_x$
go to the classes denoted by the same symbols,
the class $\omega_1$ goes to zero,
the class $\omega_3$ to $\omega_z \omega^2$,
and the class $c_4$ to $\omega^4$.
Moreover, the relation (1.5.5$'$) passes to the relation
$$
\omega_z^2 \omega^4 = \omega_z^2 \omega^4 + \gamma \omega_z^2 \omega^2 c_x
$$
whence $\gamma = 0$.
Hence the relation (1.5.5) holds in the cohomology of $G$.
This completes the proof of Theorem 1.5.

\noindent{Remark 3.7.}
Our method involving the group $\widehat G$  
or, more precisely, 
the fiber bundle (3.2),
rather than the obvious split
extension of $N$ by $K_x \times K_z$,
leads to a somewhat more economical approach
to solving the multiplicative extension
problem since then the only ambiguities
for this extension problem
come from the powers of the generator $c_x$.
The ambiguities for the multiplicative
extension problem arising from the split extension
of $N$ by $K_x \times K_z$ 
come from arbitrary monomials in the generators
$\omega_x$, $c_x$, and $\omega_z$.

\noindent
{\smc Remark 3.8.} The group $G$ may as well be written as a split extension
$G=G_x \rrtimes K_z$ and, in view of Proposition 2.1.4,
with reference to the induced filtration
of the mod 2 cohomology algebra, the associated graded algebra
is just the tensor product of 
$\roman H^*(G_x,\Bobb Z/2)$ and
$\roman H^*(K_z,\Bobb Z/2)$.
The relations (1.5.4) and (1.5.5) show that the multiplicative
extension problem is non-trivial, though.

\noindent
{\smc Remark 3.9.} The group $N$ has order $2^\rho$,
and we have made the assumption that $\rho \geq 3$.
Theorem 1.5 shows that the mod 2 cohomology algebra
of the holomorph $N \rrtimes \roman{Aut}(N)$
is then independent of $\rho$.
Given $\rho$, consider the group $N^{\sharp}= \Bobb Z\big/2^{\rho+1}$.
The obvious surjection from $N^{\sharp}$ to $N=\Bobb Z\big/2^{\rho}$
induces a surjection from $N^{\sharp} \rrtimes \roman{Aut}(N^{\sharp})$
to $N \rrtimes \roman{Aut}(N)$ but, beware,
this surjection does {\it not\/} induce an isomorphism
between the mod 2 cohomology algebras.
Indeed, under the induced morphism in cohomology,
the multiplicative generators
$c_x$, $\omega_3$, and $c_4$ go to zero.

\bigskip
\widestnumber\key{999}
\centerline{References}
\smallskip\noindent

\ref \no \differen
\by J. Huebschmann
\paper Automorphisms of group extensions and
 differentials in the Lyndon-Hochschild-Serre spectral sequence
\jour J. of Algebra
\vol 72
\yr 1981
\pages 296--334
\endref

\ref \no \modpcoho
\by J. Huebschmann
\paper The mod $p$ cohomology rings of metacyclic groups
\jour J. of Pure and Applied Algebra
\vol 60
\yr 1989
\pages 53--105
\endref

\ref \no \huebkade
\by J. Huebschmann and T. Kadeishvili
\paper Small models for chain algebras
\jour Math. Z.
\vol 207
\yr 1991
\pages 245--280
\endref

\enddocument